\documentclass[11pt]{article}
\usepackage[dvipdfmx]{graphicx}

\usepackage[
dvipdfmx,
bookmarks=true,
bookmarksnumbered=true, 
bookmarkstype=toc,
setpagesize=false,
pdftitle={Fast implementation for semidefinite programs with positive matrix completion},
pdfauthor={Makoto Yamashita},
pdfkeywords={semidefinite programs,
the primal-dual interior-point methods,
positive definite matrix completion,
multithreaded parallel computing}
]{hyperref}

\ifx\kanjiskip\undefined\else
  \usepackage{atbegshi}
  \ifx\ucs\undefined
    \ifnum 42146=\euc"A4A2
      \AtBeginShipoutFirst{}
    \else
      \AtBeginShipoutFirst{}
    \fi
  \else
    \AtBeginShipoutFirst{}
  \fi
\fi

\usepackage{amssymb}
\def\@themcountersep{}
\catcode`@=11
\renewcommand{\@begintheorem}[2]{\it \trivlist
        \item[\hskip \labelsep{\bf #1\ #2{.}}]}
\renewcommand{\@opargbegintheorem}[3]{\it \trivlist
        \item[\hskip \labelsep{\bf #1\ #2\ (#3) {:}}]}
\catcode`@=12
\newtheorem{THEO}{Theorem}[section]
\newtheorem{ALGo}[THEO]{Algorithm}
\newtheorem{CONJ}[THEO]{Conjecture}
\newtheorem{COND}[THEO]{Condition}
\newtheorem{CORO}[THEO]{Corollary}
\newtheorem{DEFI}[THEO]{Definition}
\newtheorem{EXAMP}[THEO]{Example}
\newtheorem{FACT}[THEO]{Fact}
\newtheorem{HYPO}[THEO]{Hypothesis}
\newtheorem{LEMM}[THEO]{Lemma}
\newtheorem{PROB}[THEO]{Problem}
\newtheorem{PROP}[THEO]{Proposition}
\newtheorem{REMA}[THEO]{Remark}
\newcommand{\theo}{\begin{THEO}}
\newcommand{\algo}{\begin{ALGo} \rm}
\newcommand{\cond}{\begin{COND}}
\newcommand{\conj}{\begin{CONJ}}
\newcommand{\coro}{\begin{CORO}}
\newcommand{\defi}{\begin{DEFI} \rm}
\newcommand{\examp}{\begin{EXAMP} \rm}
\newcommand{\fact}{\begin{FACT}}
\newcommand{\hypo}{\begin{HYPO} \rm}
\newcommand{\lemm}{\begin{LEMM}}
\newcommand{\prob}{\begin{PROB} \rm}
\newcommand{\prop}{\begin{PROP}}
\newcommand{\rema}{\begin{REMA} \rm}
\newcommand{\etheo}{\end{THEO}}
\newcommand{\ealgo}{\end{ALGo}}
\newcommand{\econd}{\end{COND}}
\newcommand{\econj}{\end{CONJ}}
\newcommand{\ecoro}{\end{CORO}}
\newcommand{\edefi}{\end{DEFI}}
\newcommand{\eexamp}{\end{EXAMP}}
\newcommand{\efact}{\end{FACT}}
\newcommand{\ehypo}{\end{HYPO}}
\newcommand{\elemm}{\end{LEMM}}
\newcommand{\eprob}{\end{PROB}}
\newcommand{\eprop}{\end{PROP}}
\newcommand{\erema}{\end{REMA}}

\def\0{\mbox{\bf 0}}
\def\1{\mbox{\bf 1}}
\def\2{\mbox{\bf 2}}
\def\3{\mbox{\bf 3}}
\def\4{\mbox{\bf 4}}
\def\5{\mbox{\bf 5}}
\def\6{\mbox{\bf 6}}
\def\7{\mbox{\bf 7}}
\def\8{\mbox{\bf 8}}
\def\9{\mbox{\bf 9}}

\newdimen\zhige \zhige=0pt
\def\chige#1{{\setbox\zhige\hbox{#1}\ifdim\ht\zhige=1ex\accent24 #1%
  \else\ooalign{\unhbox\zhige\crcr\hidewidth\char24\hidewidth}\fi}}

\def\d{\mbox{\boldmath $d$}}
\def\e{\mbox{\boldmath $e$}}

\def\g{\mbox{\boldmath $g$}}

\def\v{\mbox{\boldmath $v$}}
\def\w{\mbox{\boldmath $w$}}

\def\z{\mbox{\boldmath $z$}}
\def\A{\mbox{\boldmath $A$}}
\def\B{\mbox{\boldmath $B$}}

\def\D{\mbox{\boldmath $D$}}

\def\G{\mbox{\boldmath $G$}}

\def\I{\mbox{\boldmath $I$}}

\def\L{\mbox{\boldmath $L$}}
\def\M{\mbox{\boldmath $M$}}
\def\N{\mbox{\boldmath $N$}}
\def\O{\mbox{\boldmath $O$}}
\def\P{\mbox{\boldmath $P$}}

\def\U{\mbox{\boldmath $U$}}
\def\V{\mbox{\boldmath $V$}}
\def\W{\mbox{\boldmath $W$}}
\def\X{\mbox{\boldmath $X$}}
\def\Y{\mbox{\boldmath $Y$}}

\def\AC{\mbox{$\cal A$}}

\def\DC{\mbox{$\cal D$}}
\def\EC{\mbox{$\cal E$}}

\def\PC{\mbox{$\cal P$}}

\def\SC{\mbox{$\cal S$}}

\def\Real{\mbox{$\mathbb{R}$}}

\def\SMAT{\mbox{$\mathbb{S}$}}

\usepackage{enumerate}

\begin{document}

\begin{center}
 {\Large \bf 
 Fast implementation 
for semidefinite programs with positive matrix completion}
 \\
 Makoto Yamashita
 \footnote{
 Department of Mathematical and Computing Sciences,
 Tokyo Institute of Technology, 2-12-1-W8-29 Ookayama, Meguro-ku, Tokyo
 152-8552, Japan (Makoto.Yamashita@is.titech.ac.jp).
 The work of the first author was financially supported
 by the Sasakawa Scientific Research Grant from The Japan Science Society.
 },
 Kazuhide Nakata
 \footnote{
 Department of Industrial Engineering and Management,
 Tokyo Institute of Technology, 2-12-1-W9-60 Ookayama, Meguro-ku, Tokyo
 152-8552, Japan (nakata.k.ac@m.titech.ac.jp).
 The work of the second author was partially supported by
 Grant-in-Aid for Young Scientists (B) 22710136.
 }
 \\
 November 1, 2013, Revised May 26, 2014.
 \end{center}

\noindent
{\bf Abstract}

Solving semidefinite programs (SDP) in a short time is the key to
managing various mathematical optimization problems.
The matrix-completion primal-dual interior-point method (MC-PDIPM)
extracts a sparse structure of input SDP by factorizing the variable matrices. 
In this paper, we
propose a new factorization based on the inverse of the variable matrix
to enhance the performance of MC-PDIPM.  We also use
multithreaded parallel computing to deal with the major bottlenecks in 
MC-PDIPM.  Numerical results show that the new factorization and 
multithreaded computing reduce the computation time for 
SDPs that have structural sparsity.

\vspace{0.5cm}

\noindent
{\bf AMS Classification} \\ 
90~Operations~research,~mathematical~programming,
90C22~Semidefinite~programming, 
90C51~Interior-point~methods,
97N80~Mathematical~software,~computer~programs.

\section{Introduction}

Semidefinite programs (SDP) have become one of main topics of
mathematical optimization, because of  its wide range application 
from combinatorial optimization~\cite{GOEMANS95} to 
quantum chemistry~\cite{FUKUDA07,NAKATA01}
and sensor network localization~\cite{BISWAS04}.
A survey of its many applications can be found in Todd's paper~\cite{TODD01},
and the range is still expanding.
Moreover, there is no doubt that solving SDPs in a short time 
is the key to managing such applications.
The primal-dual interior-point method
(PDIPM)~\cite{ALIZADEH98, HELMBERG96, KOJIMA94, MONTEIRO97, NESTEROV98}
is often employed since it can solve SDPs in a polynomial time, and 
many solvers are based on it, for example,
SDPA~\cite{SDP-HANDBOOK}, CSDP~\cite{BORCHERS99}, 
SeDuMi~\cite{STURM99}, and SDPT3~\cite{TODD99}.
A recent paper \cite{YAMASHITA12} reports that  integration with parallel computing
enables one to solve 
large-scale SDPs arising in practical applications.
% A recent paper in 2012~\cite{FUJISAWA12} reports that
% we, with the help of super-computers, are really ready
% to solve an extremely large-scale SDP whose size is far beyond the SDPs 
% we challenged a few years ago.

A major difficulty with PDIPM is that 
the primal variable matrix $\X$ must be handled as a fully dense matrix
even when
all the input data matrices $\A_0,\ldots,\A_m$ are considerably sparse.
The standard form in this paper is the primal-dual pair;
\begin{eqnarray*}
 (\PC) & & \begin{array}{lcl}
  \min &:& \A_0 \bullet \X \\
	    \mbox{subject to} &:& \A_k \bullet \X = b_k \quad
	     (k=1,\ldots,m) \\
	    &:& \X \succeq \O 
	     \end{array} \\
 (\DC) & & \begin{array}{lcl}
  \max &:& \sum_{k=1}^m b_k z_k \\
	    \mbox{subject to} &:& \sum_{k=1}^m \A_k z_k + \Y = \A_0 \\
	    &:& \Y \succeq \O 
	     \end{array} 
 \end{eqnarray*}
Let $\SMAT^n$ be the space of $n \times n$ symmetric matrices.
The symbol $\X \succeq \O (\X \succ \O)$ indicates that  $\X \in \SMAT^n$ is 
a positive semidefinite (definite) matrix.
The notation $\U \bullet \V$ is the inner-product between
$\U,\V \in \SMAT^n$ defined by
$\U \bullet \V = \sum_{i=1}^n \sum_{j=1}^n U_{ij} V_{ij}$.
The input data are $\A_0,\A_1,\ldots,\A_m \in \SMAT^n$ and
$b_1,\ldots,b_m \in \Real$.
The variable in the primal problem $(\PC)$ is
$\X \in \SMAT^n$, while the variable in the dual problem $(\DC)$ 
is $\Y \in \SMAT^n$ and $\z \in \Real^m$.

The sparsity of the input matrices directly affects 
the dual matrix $\Y = \A_0 - \sum_{k=1}^m \A_k z_k$.
More precisely, $Y_{ij}$ can be nonzero only when 
the {\it aggregate sparsity pattern} defined by
$\AC = \{(i,j) : 1 \le i \le n, 1 \le j \le n,
[\A_k]_{ij} \ne 0 \ \mbox{for some} \ k=0,\ldots,m\}$ covers $(i,j)$.
Here, $[\A_k]_{ij}$ is the $(i,j)$th element of $\A_k$.
Examples of aggregate sparsity patterns are illustrated 
in Figures~\ref{fig:max-clique} and \ref{fig:spin-glass}; 
these patterns arise from the SDPs we solved in the numerical 
experiments.
%It should be emphasized that we can ignore $Y_{ij}$
%when $(i,j) \notin \AC$, since such $Y_{ij}$ should be zero.
On the other hand, 
all the elements of $\X$  in the primal problem $(\PC)$
must be stored in memory in order to check the constraints $\X \succeq \O$.
The matrix-completion primal-dual interior-point method (MC-PDIPM)
proposed in \cite{FUKUDA00, NAKATA03} enables the PDIPM
to be executed by factoring $\X$ into the form,
\begin{eqnarray}
\X = \L_1^T \L_2^T \cdots \L_{\ell-1}^T \D \L_{\ell-1} \cdots \L_2 \L_1
\label{eq:original-decomp}
\end{eqnarray}
where $\D$ is a diagonal-block positive semidefinite matrix
and $\L_1, \L_2,\ldots, \L_{\ell-1}$ are lower triangular matrices.
A remarkable feature of this factorization is that 
$\D$ and $\L_1, \L_2,\ldots, \L_{\ell-1}$ inherit the sparsity of $\AC$.
When $\AC$ is considerably sparse, these matrices are also sparse;
hence, MC-PDIPM has 
considerable advantages compared with handling the fully dense matrix $\X$.
It was first implemented in a solver called
SDPA-C (SemiDefinite Programming Algorithm with Completion), 
a variant of SDPA~\cite{SDP-HANDBOOK},
and as reported in \cite{FUKUDA00, NAKATA03},
it significantly reduces computation costs of solving
SDPs with structural sparsity.
% compared to the standard PDIPM.

The main objective of this paper is to accelerate MC-PDIPM.
The principal bottleneck is 
the repeated computation of the form $\X \v$ for $\v \in \Real^n$. 
The original factorization (\ref{eq:original-decomp}) can be summarized as 
$\X = \L^{\ T} \D \L$ with a lower triangular
matrix $\L = \L_{\ell-1}\cdots\L_2\L_1$.
Instead of this factorization,
we introduce the Cholesky factorization of the inverse of $\X$;
$\X^{-1} = \widehat{\L}\widehat{\L}^{T}$
and show that the lower triangular matrix $\widehat{\L}$ directly
inherits the sparsity from $\AC$. % in a better way than $\D$ and $\L$.
Another obstacle of (\ref{eq:original-decomp}) is that 
the presence of $\D$ means that it is not 
a standard form of  Cholesky factorization, and this
prevents us from using software packages that are available for sparse Cholesky
factorization, such as CHOLMOD~\cite{CHOLMOD} and MUMPS~\cite{MUMPS-A}. 
However, the removal of $\D$ by using $\X^{-1} = \widehat{\L}\widehat{\L}^T$
would enable us to naturally integrate these packages into MC-PDIPM framework;
in so doing, we can obtain the results of $\X \v$ in a more effective way 
and shrink the computation time of MC-PDIPM.

In this paper, we also introduce multithreaded parallel computing to this new factorization.
Most processors on modern PCs have multiple cores,
and we can process some tasks simultaneously on different cores.
A parallel computation of MC-PDIPM on multiple PCs connected by a
local area network was already discussed in \cite{NAKATA06}.
In this paper, we employ different parallel schemes
for multithreading on a single PC, because the differences between the memory accesses
of  parallel computing
with the Message Passing Interface (MPI) protocol
on multiple PCs and those of multithreading on a single PC strongly affects
the performance of parallel computing.
In addition, to enhance the performance of multithreading,
we control the number of threads involved in our parallel schemes.

On the basis of  the existing version SDPA-C~6.2.1, we implemented a new version,
SDPA-C~7.3.8
(The version numbers reflect the versions of SDPA that SDPA-C branches from).
We conducted numerical experiments, showing that the new SDPA-C~7.3.8
successfully reduces the computation time
because of the effectiveness of $\X^{-1} = \widehat{\L}\widehat{\L}^T$.
We also show that the multithreaded computation further
expands the difference in computation time between SDPA-C~6.2.1 and
7.3.8.

This paper is organized as follows.
In Section~\ref{sec:pre}, we introduce two preliminary concepts,{\it i.e.}, 
the positive matrix completion and PDIPM.
Section~\ref{sec:new} is the main part of this paper that
describes the new implementation in detail.
Section~\ref{sec:results} presents numerical results showing
its performance.
In Section~\ref{sec:conclusions}, we summarize this paper
and discuss future directions.

Throughout this paper, 
we will use $|S|$ to denote the number of elements of the set $S$.
For a matrix $\X$ and two sets $S,T \subset \{1,\ldots,n\}$,
we use the notation $\X_{ST}$ to denote the sub-matrix of $\X$ that
collects the elements of $X_{ij}$ with $i \in S$ and $j \in T$; for
example, $\X_{\{2,6\},\{3,4\}} = \left(\begin{array}{cc}
				  X_{23} & X_{24} \\ X_{63} & X_{64} 
				   \end{array}\right)$.

\section{Preliminaries}\label{sec:pre}

Here, we briefly describe the basic concepts of  positive matrix completion 
and PDIPM.
For more details on the two and their relation,
please refer to \cite{FUKUDA00, NAKATA03} and references therein.

\subsection{Positive Matrix Completion}\label{sec:pmc}

Positive matrix completion is closely related to the Cholesky
factorization of the variable matrices $\X$ and $\Y$ in 
the context of the PDIPM framework.
When $\Y = \A_0 - \sum_{k=1}^m \A_k z_k$ is positive definite, 
we can apply Cholesky factorization to obtain 
a lower triangular matrix $\N$ such that $\Y = \N \N^T$. 
However, this factorization generates nonzero elements
out of the aggregate sparsity pattern $\AC$, and this phenomenon is called fill-in. 
Although $\AC$ is not
enough to cover all the nonzeros in $\N$,
it is known that
we can prepare a set of appropriate subsets
$C_1,\ldots,C_{\ell} \subset \{1,2,\ldots,n\}$
so that the set 
$ \EC = \cup_{r=1}^{\ell} (C_r \times C_r)$
covers the nonzero positions of $\AC$ and the fill-in.
These subsets $C_1,\ldots,C_{\ell} \subset \{1,2,\ldots,n\}$ 
are called cliques, in relation with graph theory,
and they are obtained
in three steps; we permute the rows/columns of $\Y$ with an
appropriate order like approximation minimum ordering and
generate a chordal graph from $\AC$. Then, we extract the maximal
cliques there as $C_1,\ldots,C_{\ell}$.
The set $\EC$ is called
the {\it extended sparsity pattern}.
Throughout this paper, we will assume that $\EC$ is considerably sparse;
$|\AC|$ and $|\EC|$
are much less than the fully dense case
$n^2$, for instance, $|\AC| \le |\EC| < 10^{-2} \times n^2$
for large $n$.
In addition, we will assume for simplicity
that $C_1,\ldots,C_{\ell} \subset \{1,2,\ldots,n\}$ 
are sorted in an appropriate order which satisfies a nice property,
called the running intersection property in \cite{FUKUDA00}.
Such an order can be easily derived from the chordal graph.

Grone {\it et al.}~\cite{GRONE84} proved that if a given matrix $\overline{\X}$
satisfies the positive-definite conditions on all the sub-matrices
induced by the cliques $C_1,C_2,\ldots,C_{\ell}$; that is,
$\overline{\X}$ satisfies
$\overline{\X}_{C_1C_1} \succ \O, \overline{\X}_{C_2C_2} \succ \O, \ldots,
\overline{\X}_{C_{\ell}C_{\ell}} \succ \O$, then
$\overline{\X}$ can be completed to
$\overline{\overline{\X}}$ such that 
$\overline{\overline{\X}}_{C_rC_r} = \overline{\X}_{C_rC_r}$ for
$r=1,\ldots, \ell$
and the entire matrix $\overline{\overline{\X}}$ is positive definite.
Furthermore, it was shown in \cite{FUKUDA00} that 
the explicit formula (\ref{eq:sparse_fact}) below completes $\overline{\X}$ 
to the max-determinant completion $\widehat{\X}$, which satisfies
\begin{eqnarray*}
 \det(\widehat{\X}) = \max \{ \det(\overline{\overline{\X}}) : 
\overline{\overline{\X}}_{C_rC_r} =
\overline{\X}_{C_rC_r}\ \mbox{for} \ r=1,\ldots,\ell,
\quad \overline{\overline{\X}} \succ \O\}.
\end{eqnarray*}
The sparse factorization of $\widehat{\X}$ from $\overline{\X}$ is given by
\begin{eqnarray}
\widehat{\X} = \L_1^T \L_2^T \ldots \L_{\ell-1}^T \D \L_{\ell-1} \ldots \L_2 \L_1,
\label{eq:sparse_fact}
\end{eqnarray}
where $\L_1,\L_2,\ldots,\L_{\ell-1}$ are the triangular lower matrices of the form,
\begin{eqnarray}
 \left[\L_r\right]_{ij} = \left\{
		\begin{array}{ll}
		 1 & (i=j) \\
		 \left[\overline{\X}^{-1}_{U_rU_r}
		  \overline{\X}_{U_rS_r}\right]_{ij} & (i \in U_r, j \in S_r) \\
		 0 & (otherwise)
		\end{array}\right.
 \label{eq:orig_L}
\end{eqnarray}
and $\D$ is the diagonal-block matrix 
\begin{eqnarray}
 \D = \left(\begin{array}{cccc}
       \D_{S_1S_1} & & & \\
	     & \D_{S_2S_2} & & \\
	     & & \ddots & \\
	     & & & \D_{S_{\ell}S_{\ell}} 
  \end{array}\right)
\end{eqnarray} \label{eq:orig_D}
with
\begin{eqnarray*}
 S_r &=& C_r \backslash (C_{r+1} \cup C_{r+2} \cup \cdots \cup C_{\ell})
  \quad (r=1,2,\ldots,\ell) \\
 U_r &=& C_r \cap (C_{r+1} \cup C_{r+2} \cup \cdots \cup C_{\ell})
  \quad (r=1,2,\ldots,\ell)
\end{eqnarray*}
and
\begin{eqnarray*}
 \D_{S_rS_r} = \left\{\begin{array}{ll}
		\overline{\X}_{S_rS_r}
		 - \overline{\X}_{S_rU_r}\overline{\X}_{U_rU_r}^{-1}
		 \overline{\X}_{U_rS_r} & (r=1,2,\ldots,\ell-1) \\
		       \overline{\X}_{S_{\ell}S_{\ell}} & (r=\ell) 
\end{array}\right..
\end{eqnarray*}

It can be shown that the triangular lower matrix $\L$ defined by
$\L = \L_{\ell-1} \ldots \L_2 \L_1$ is usually fully dense,
thereby destroying the structural sparsity of $\EC$.
%The completed matrix $\widehat{\X}$ also has nonzeros out of $\EC$,
%and should be handled as a fully dense matrix.
Therefore, when we compute
$\w = \widehat{\X} \v = \L^T \D \L \v$
for some vector $\v \in \Real^m$, constructing a fully dense $\L$ is not 
efficient. 
We should note that the inverse $\L^{-1}$ can maintain the sparsity of $\EC$, that is,
$[\L^{-1}]_{ij} = 0$ if $(i,j) \notin \EC$, and $\L^{-1}$ is a lower
triangular matrix~\cite{NAKATA03}.
Hence, the two equations 
$\L^{-1} \w_1 = \v $ and $\L^{-T}\w = \D \w_1$
can be solved with forward/backward substitutions
by exploiting the structure of $\EC$, 
and they can be used to compute $\w$ much faster.
In addition, we  do not need to compose a fully dense
$\widehat{\X}$ via multiplication of $\L^T \D \L$.
%than composing fully dense
%$\widehat{\X}$ via the multiplication of $\L^T \D \L$.
This idea saves on the computation cost of PDIPM, as
discussed in the next subsection.

\subsection{Primal-Dual Interior-Point Method}

This subsection briefly describes the primal-dual interior-point method (PDIPM)
and the modification of its computation formula
by using the positive matrix completion method.

The basic framework of the PDIPM can be summarized as follow.

\vspace{0.5cm}
\noindent{\bf Basic framework of the primal-dual interior-point
method}
\begin{enumerate}
 \item[Step 0] Choose an initial point $(\X, \Y, \z)$ such that
	      $\X \succ \O$ and $\Y \succ \O$.
	      Choose parameters $\beta$ and $\gamma$ from $0 < \beta < 1$
	      and $0 < \gamma < 1$.
 \item[Step 1] If $(\X, \Y, \z)$ satisfies a stopping criterion,
	      output $(\X, \Y, \z)$ as a solution and terminate.
 \item[Step 2] Compute a search direction $(d\X, d\Y, d\z)$
	      based on the modified Newton method.
 \item[Step 3] Compute the maximum step length $\alpha_p$ and $\alpha_d$ such that
	      \begin{eqnarray}
	       \alpha_p &=& \max\{\alpha \in (0,1] :
		\X + \alpha d\X \succeq \O\} \label{eq:primal_length} \\
	       \alpha_d &=& \max\{\alpha \in (0,1] :
		\Y + \alpha d\Y \succeq \O\}. \nonumber
	      \end{eqnarray}
 \item[Step 4] Update $(\X, \Y, \z)$ with
	      $(\X + \gamma \alpha_p d\X, \Y + \gamma \alpha_d d\Y, \z + \gamma \alpha_d d\z)$.
	      Go to Step 1.
\end{enumerate}

\vspace{0.5cm}
The chief computation in the above framework is usually
that of the search direction $(d\X, d\Y, d\z)$, as pointed out in \cite{YAMASHITA02}.
If we employ the HKM direction~\cite{HELMBERG96,KOJIMA94,MONTEIRO97},
the search direction can be obtained with the following system;
\begin{eqnarray}
 \B \d\z &=& \g \label{eq:SCE}\\
 d\Y &=& \G - \sum_{k=1}^m \A_k dz_k \nonumber \\
 \widehat{d\X} &=& \beta \mu \Y^{-1} - \X - \X d\Y \Y^{-1} \label{eq:dX}, d\X = (\widehat{d\X} + \widehat{d\X}^T)/2
\end{eqnarray}
where
\begin{eqnarray}
 B_{ij} &=& (\X \A_i \Y^{-1}) \bullet \A_j \quad
  (i=1,\ldots,m,\ j=1,\ldots,m) \label{eq:SCM} \\
 g_k &=& \A_k \bullet (\beta \mu \Y^{-1} - \X - \X \G \Y^{-1}) \quad (k=1,\ldots,m) \nonumber
\end{eqnarray}
with $\mu =\frac{\X \bullet \Y}{n}$, 
$\G = \A_0 - \sum_{k=1}^m \A_k z_k$.
The linear system (\ref{eq:SCE}) is often called
the {\it Schur complement equation} (SCE), and its coefficient matrix
determined from (\ref{eq:SCM}) is called the {\it Schur complement matrix} (SCM).
We first solve  SCE (\ref{eq:SCE}) to obtain $d\z$ and then
compute $d\Y$ and $d\X$.

The matrix completion (\ref{eq:original-decomp})
enables us to replace the fully dense
matrices $\widehat{\X}$ and $\Y^{-1}$ with their sparse versions
in the above computation.
From the properties of the inner product,
the change from $\X$ to  $\widehat{\X}$ in formula (\ref{eq:SCM}) 
does not affect $B_{ij}$; therefore,
its formula can be transformed into
\begin{eqnarray}
 B_{ij} &=& (\widehat{\X} \A_i \Y^{-1}) \bullet \A_j \nonumber \\
  &=& \sum_{k=1}^m (\L^T \D \L \e_k)^T
   \A_i (\N^{-T}\N^{-1}[\A_j]_{*k}) \label{eq:newSCM}
\end{eqnarray}
where $\e_k$ and $[\A_j]_{*k}$ are the $k$th columns of $\I$ and $\A_j$,
respectively. 

We also modify the computation of the primal search direction $d\X$
by evaluating its auxiliary matrix $\widehat{d\X}$ in a column-wise manner
\begin{eqnarray}
 [\widehat{d\X}]_{*k} &=&
  \beta \mu \Y^{-1} \e_k - \widehat{\X}\e_k - \widehat{\X} d\Y
  \Y^{-1}\e_k  \nonumber \\
&=&  \beta \mu \N^{-T}\N^{-1} \e_k
   - \L^T \D \L \e_k
   - \L^T \D \L d\Y
   \N^{-T}\N^{-1} \e_k \label{eq:newdX}.
 \end{eqnarray}
As pointed out in Section \ref{sec:pmc}, 
by solving 
the linear equations that involve the sparse matrices $\L^{-1}$
and $\N$, % with their forward/backward substitutions.
we can avoid the fully dense matrices $\widehat{\X}$ and $\Y^{-1}$ 
in (\ref{eq:newSCM}) and (\ref{eq:newdX}).
The computation of the step length $\alpha_p$ in 
(\ref{eq:primal_length}) can also be decomposed into the sub-matrices
\begin{eqnarray}
 \widehat{\alpha}_p = \min_{r=1,2,\ldots,\ell}
  \max\{\alpha \in (0,1] : \overline{\X}_{C_rC_r} + \alpha d\X_{C_rC_r} \succeq \O\}
\label{eq:primal_length2},
\end{eqnarray}
so that
$\overline{\X}_{C_rC_r} + \widehat{\alpha}_p d\X_{C_rC_r}$ 
is positive definite for $r=1,\ldots,\ell$,
and we can complete these sub-matrices to the positive definite
matrix $\widehat{\X}$.

The numerical results in \cite{NAKATA03} indicated that 
removal of the fully dense matrices $\widehat{\X}$ and $\Y^{-1}$
makes the MC-PDIPM run more effectively than the standard PDIPM
({\it i.e.}, a PDIPM which does not use the positive matrix completion method)
for some types of SDP that have the structural sparsity in $\EC$.

\section{Fast implementation of the matrix-completion primal-dual
 interior-point method}\label{sec:new}

MC-PDIPM was first implemented in the solver SDPA-C~5~\cite{NAKATA03}.
Along with the update of SDPA based on the standard PDIPM to version 6,
SDPA-C was also updated to SDPA-C~6.
SDPA-C~6 utilizes the BLAS (Basic Linear Algebra Subprograms)
library~\cite{LAWSON79}
to accelerate the linear algebra computation involved in MC-PDIPM.

The new SDPA-C, version 7.3.8, described in this paper 
further reduces the computation time from that of version 6.2.1.
In this section, we describe the new features of SDPA-C~7.3.8;
the improvements in the factorization of $\widehat{\X}$ are in Section~3.1, and 
the multithreaded parallel computing for SCM $\B$ and
the primal auxiliary direction $\widehat{d\X}$ are in Section~3.2.
In what follows, we will abbreviate SDPA-C~6.2.1 and SDPA-C~7.3.8
to SDPA-C~6 and SDPA-C~7, respectively.

\subsection{New Factorization of the Completed Matrix}\label{sec:decomp}

The factorization of $\widehat{\X}$ into
$\widehat{\X} = \L^T \D \L$ is not a standard Cholesky
factorization due to the diagonal-block matrix $\D$; hence,
we could not employ software packages for the sparse Cholesky
factorization.
The completed matrix $\widehat{\X}$ is usually fully dense,
while the sparsity of its inverse $\widehat{\X}^{-1}$ inherits 
the structure of $\EC$, {\it i.e.}, $[\widehat{\X}^{-1}]_{ij} = 0$
for $(i,j) \notin \EC$.
Therefore, we will focus on $\widehat{\X}^{-1}$ rather than $\widehat{\X}$
and introduce a new factorization of the form
$\widehat{\X}^{-1} = \widehat{\L}\widehat{\L}^T$
with the lower-triangular matrix $\widehat{\L}$.
We want to emphasize here that $\widehat{\L}$ also inherits
the structure of $\EC$. 
In this subsection, we show that we can obtain
the factorized matrix $\widehat{\L}$ from $\overline{\X}$
in an efficient way by using
the structure of $S_r$ and $C_r$ $(r=1,\ldots,\ell)$.

Algorithm~1 is used to obtain $\widehat{\L}$.
The input is $\overline{\X}$, and since $\overline{\X}$
is going to be completed to a positive definite matrix $\widehat{\X}$,
we suppose that $\overline{\X}_{C_rC_r} \succ \O$ $(r=1,\ldots,\ell)$.
The validity of the algorithm will be discussed later.

\vspace{0.5cm}
\noindent{\bf Algorithm 1: Efficient algorithm to obtain
the Cholesky factorization of the inverse of the completed matrix}
\begin{enumerate}
 \item[Step 1] Initialize the memory space for $\widehat{\L}$ by
	      $\EC = \cup_{r=1}^{\ell} (C_r \times C_r)$.
 \item[Step 2] For $r=1,\ldots,\ell$, apply  Cholesky factorization
	      to $\overline{\X}^{-1}_{C_rC_r}$ to
	      obtain the lower triangular matrix $\L_r$ that satisfies
	      $\overline{\X}^{-1}_{C_rC_r} = \L_r \L_r^T$.
            We take the following steps to avoid computing $\overline{\X}^{-1}_{C_rC_r}$,
	      \begin{enumerate}
	       \item[Step 2-1]
			    Let $\P_r$ be the permutation matrix
			    of dimension
			    $|C_r| \times |C_r|$ with
			    \begin{eqnarray*}
			     \left\{\begin{array}{ll}
			     \left[\P_r\right]_{ij} = 1
			      & \mbox{if} \quad i+j=|C_r|+1 \\
			     \left[\P_r\right]_{ij} = 0  & \mbox{otherwise}
			      \end{array}\right.
			    \end{eqnarray*}
			    so that $\P_r \overline{\X}_{C_rC_r} \P_r^T$
			    has the inverse row/column order of
			    $\overline{\X}_{C_rC_r}$.
	       \item[Step 2-2]
			    Apply Cholesky factorization
			    to $\P_r \overline{\X}_{C_rC_r} \P_r^T$
			    to obtain a lower triangular matrix
			    $\M_r$ that satisfies
			    $\P_r \overline{\X}_{C_rC_r} \P_r^T = \M_r \M_r^T$.
	       \item[Step 2-3]
			    Let $\L_r$ be $\P_r \M_r^{-T} \P_r^T$.
	      \end{enumerate}
 \item[Step 3] For $r=1,\ldots,\ell$, put
	      the first $|S_r|$ columns of
	      $\L_r$ in the memory space of $\widehat{\L}_{C_rS_r}$.
\end{enumerate}

Algorithm~1 requires neither a fully dense $\widehat{\X}$ nor its
inverse $\widehat{\X}^{-1}$. In addition, since most of the computation is devoted to
the Cholesky factorization of $\P_r \overline{\X}_{C_rC_r} \P_r^T$, we can
expect there will be a considerable reduction in computation time
when the extended sparsity pattern $\EC$ is decomposed into
small $C_1,C_2,\ldots,C_{\ell}$.
Furthermore, Algorithm~1 assures that all nonzero elements
of $\widehat{\L}$ appear only in $\EC$.

\noindent{\it Validity of Algorithm~1}:

We will prove the validity of Algorithm~1 on the basis of
Lemma~2.6 of \cite{FUKUDA00}.
For simplicity, we will focus on the first clique $C_1$
and wrap up the other cliques into $C'_2 := \cup_{r=2}^{\ell} C_r$.
Because of the running intersection property, the other cliques $C_2,\ldots,C_r$
can be handled in the same way
by induction on the number of cliques.
Thanks to this property as well, we can suppose that $i<j$ for
$i \in C_1$ and $j \in C'_2 \backslash C_1$.
For $\overline{\X} \in \SMAT^n$, we decompose $\{1,2,\ldots,n\}$
into three sets 
$S = C_1 \backslash C'_2,
U = C_1 \cap C'_2$ and $T = C'_2 \backslash C_1$.
Note that the extended sparsity pattern $\EC$ of $\overline{\X}$ is
covered by
$\left(\left(S \cup U\right) \times \left(S \cup U\right) \right)
\cup
\left(\left(U \cup T\right) \times \left(U \cup T\right) \right)$.
Hence, the situation is one where $\overline{\X}$ is of the form,
\begin{eqnarray*}
 \overline{\X} = \left(\begin{array}{ccc}
		  \overline{\X}_{SS} & \overline{\X}_{SU} &
		       ? \\
		  \overline{\X}_{US} & \overline{\X}_{UU} &
		       \overline{\X}_{UT} \\
		  ? & \overline{\X}_{TU} &
		       \overline{\X}_{TT} \\
		       \end{array}\right)
\end{eqnarray*}
with unknown elements $?$ in the position 
$\left(S \times T\right) \cup \left(T \times S\right)$,
and the sub-matrices induced by the cliques $C_1, C'_2$
are positive definite,
\begin{eqnarray*}
 \overline{\X}_{C_1C_1} = \left(\begin{array}{cc}
			   \overline{\X}_{SS} & \overline{\X}_{SU} \\
				 \overline{\X}_{US} & \overline{\X}_{UU} 
				\end{array}\right) \succ \O,
 \quad 
 \overline{\X}_{C'_2C'_2} = \left(\begin{array}{cc}
			   \overline{\X}_{UU} & \overline{\X}_{UT} \\
				 \overline{\X}_{TU} & \overline{\X}_{TT} 
				\end{array}\right) \succ \O.
\end{eqnarray*}
Note that $C_1 = S \cup U$ and $C'_2 = U \cup T$ in this situation.
Lemma~2.6 of \cite{FUKUDA00} claims that $\overline{\X}$ can be
completed to the max-determinant positive definite matrix
$\widehat{\X}$,
\begin{eqnarray*}
 \widehat{\X} = \left(\begin{array}{ccc}
		  \overline{\X}_{SS} & \overline{\X}_{SU} &
		       \overline{\X}_{SU}\overline{\X}^{-1}_{UU}
		       \overline{\X}_{UT} \\
		  \overline{\X}_{US} & \overline{\X}_{UU} &
		       \overline{\X}_{UT} \\
		       \overline{\X}_{TU}\overline{\X}^{-1}_{UU}
		       \overline{\X}_{US}  & \overline{\X}_{TU} &
		       \overline{\X}_{TT} \\
		       \end{array}\right) \succ \O.
\end{eqnarray*}
Hence, proving the validity of Algorithm~1 reduces to proving that
$\widehat{\X} = \widehat{\L}^{-T} \widehat{\L}^{-1}$.
In Step~2, the inverses of the positive definite 
sub-matrices are factorized into the lower triangular matrices by using Cholesky
factorization, as follow;
\begin{eqnarray*}
 \left(\begin{array}{cc}
  \overline{\X}_{SS} & \overline{\X}_{SU} \\
  \overline{\X}_{US} & \overline{\X}_{UU}
   \end{array}\right)^{-1}
 &=&
 \left(\begin{array}{cc}
  \M_{SS} &  \\
  \M_{US} & \M_{UU}
   \end{array}\right)
 \left(\begin{array}{cc}
  \M_{SS}^T & \M_{US}^T \\
   & \M_{UU}^T
   \end{array}\right) \\
 \left(\begin{array}{cc}
  \overline{\X}_{UU} & \overline{\X}_{UT} \\
  \overline{\X}_{TU} & \overline{\X}_{TT}
   \end{array}\right)^{-1}
 &=&
 \left(\begin{array}{cc}
  \N_{UU} &  \\
  \N_{TU} & \N_{TT}
   \end{array}\right)
 \left(\begin{array}{cc}
  \N_{UU}^T & \N_{TU}^T \\
   & \N_{TT}^T
   \end{array}\right).
 \end{eqnarray*}
%Since the matrcies in the right-hand side are lower-triangular,
%their inverse matrices are also lower-triangular.
Since the matrices on the left-hand side are positive definite,
we can take the inverses of components on the right-hand side,
{\it e.g.}, $\M_{SS}^{-1}$.
By comparing the elements of both sides, we obtain
\begin{eqnarray}
 \left\{\begin{array}{ccl}
 \overline{\X}_{SS} &=& \M_{SS}^{-T} \M_{SS}^{-1}
  + \M_{SS}^{-T} \M_{US}^T \M_{UU}^{-T} \M_{UU}^{-1} \M_{US} \M_{SS}^{-1} \\
 \overline{\X}_{SU} &=& -\M_{SS}^{-T}\M_{US}^T\M_{UU}^{-T}\M_{UU}^{-1} \\
 \overline{\X}_{UU} &=& \M_{UU}^{-T}\M_{UU}^{-1}
  = \N_{UU}^{-T} \N_{UU}^{-1}
  + \N_{UU}^{-T} \N_{TU}^T \N_{TT}^{-T} \N_{TT}^{-1} \N_{TU} \N_{UU}^{-1} \\
 \overline{\X}_{UT} &=& -\N_{UU}^{-T}\N_{TU}^T\N_{TT}^{-T}\N_{TT}^{-1} \\
 \overline{\X}_{TT} &=& \N_{TT}^{-T}\N_{TT}^{-1}.
\end{array}\right. \label{eq:elements-wise}
\end{eqnarray}

In Step~3, the elements of the above factorized matrices are located 
in $\widehat{\L}$ as follows:
\begin{eqnarray}
 \widehat{\L} = \left(\begin{array}{ccc}
		 \M_{SS} & & \\
		       \M_{US} & \N_{UU} & \\
		       & \N_{TU} & \N_{TT}
			\end{array}\right), \label{eq:widehatL}
\end{eqnarray}
and since $\widehat{\L}$ is lower triangular, its inverse can be explicitly
obtained as 
\begin{eqnarray*}
 \widehat{\L}^{-1} = \left(\begin{array}{ccc}
		      \M_{SS}^{-1} & & \\
			    -\N_{UU}^{-1}\M_{US}\M_{SS}^{-1}
			     & \N_{UU}^{-1} \\
			    \N_{TT}^{-1}\N_{TU}\N_{UU}^{-1}\M_{US}\M_{SS}^{-1}
			     & -\N_{TT}^{-1}\N_{TU}\N_{UU}^{-1}
			     & \N_{TT}^{-1}
			     \end{array}\right).
\end{eqnarray*}
The product
$\widehat{\L}^{-T}\widehat{\L}^{-1}$
and  (\ref{eq:elements-wise})
lead to $\widehat{\X} = \widehat{\L}^{-T}\widehat{\L}^{-1}$.
This completes the proof of validity of Algorithm~1.
$\hfill\Box$ % This box indicates the end of proof.

As a result of  this new factorization, the evaluation formula of
SCM~(\ref{eq:newSCM})
and the primal auxiliary matrix~(\ref{eq:newdX}) can be replaced with
 efficient ones, {\it i.e.}, 
\begin{eqnarray}
 B_{ij} 
  &=& \sum_{k=1}^m (\widehat{\L}^{-T}\widehat{\L}^{-1} \e_k)^T
   \A_i (\N^{-T}\N^{-1}[\A_j]_{*k}) \label{eq:newSCM2} \\
 \left[\widehat{d\X}\right]_{*k} 
&=&  \beta \mu \N^{-T}\N^{-1} \e_k
   - \widehat{\L}^{-T}\widehat{\L}^{-1} \e_k
   - \widehat{\L}^{-T}\widehat{\L}^{-1} d\Y
   \N^{-T}\N^{-1} \e_k. \label{eq:newdX2}
\end{eqnarray}

Note that once we have the factorization 
$\widehat{\X}^{-1} = \widehat{\L} \widehat{\L}^T$, we can use the SDPA-C~6
routine to compute $\widehat{\X} = \L\D\L^T$ by (\ref{eq:sparse_fact}) 
and get $\widehat{\L} = \L^{-T} \D^{-1/2}$,
where the matrix $\D^{1/2}$ is such that $\D = \D^{1/2}\D^{1/2}$.
This computation, however, requires one more step to obtain $\D^{-1/2}$ as well as
memory for $\L^{-T}, \D$ and $\D^{-1/2}$; hence, 
direct acquisition of $\widehat{\L}$ by Algorithm 1 
is more efficient than using the formula, $\widehat{\L} = \L^{-T} \D^{-1/2}$.

To implement this new factorization based on Algorithm~1 in SDPA-C~7,
five types of the computation related to Cholesky factorization  should be
employed;
\begin{enumerate}[(i)]
\item dense Cholesky factorization for 
the SCM $\B$ and its forward/backward substitution
for  SCE (\ref{eq:SCE}) if $\B$ is fully dense
\item sparse Cholesky factorization for the
SCM $\B$ and its forward/backward substitution
for SCE (\ref{eq:SCE}) if $\B$ is considerable sparse
\item dense Cholesky factorization for the sub-matrices
$\overline{\X}_{C_1C_1},\ldots,\overline{\X}_{C_{\ell}C_{\ell}}$
in Step~2 of Algorithm~1
\item  forward/backward substitution of $\widehat{\L}$ 
to solve the linear systems of the form $\widehat{\L}\widehat{\L}^T \w = \v$
\item sparse Cholesky factorization for the dual variable matrix
$\Y = \N\N^T$ and 
forward/backward substitution of $\N$ 
to solve the linear systems of form $\N\N^T \w = \v$.
\end{enumerate}
For types (i) and (ii), the sparsity of the SCM $\B$ 
heavily depends on the types of application that generates the input SDP,
as pointed out in \cite{YAMASHITA12}.
For example, SDPs arising from quantum chemistry~\cite{FUKUDA07,NAKATA01}
have fully dense SCMs; in contrast, the density of the SCMs of SDPs
arising from sensor network localization problems~\cite{KIM12}
is often less than $1\%$.
Hence, we should appropriately choose software packages for either fully dense or
sparse SCMs. We employed the dense Cholesky factorization routine
of LAPACK~\cite{LAPACK} for (i) and the sparse Cholesky factorization routine
of MUMPs~\cite{MUMPS-A} for (ii); 
one of these routines can be selected according to
the criteria proposed in \cite{YAMASHITA12} that uses information from
the input SDP.
The LAPACK routine is also applied to type (iii).

For types (iv) and (v), 
we chose CHOLMOD~\cite{CHOLMOD} rather than MUMPS~\cite{MUMPS-A},
since we must access the internal data structure of the software package 
in order to locate $\widehat{\L}_r$ in the appropriate space of $\widehat{\L}$
in Step~3 of Algorithm~1. CHOLMOD internally uses a super-nodal
Cholesky factorization, and we noticed that the row sets and the column sets 
of super-nodes computed in CHOLMOD have the running intersection
property; hence,
the row sets and columns sets can be used as
$C_1,C_2,\ldots,C_{\ell}$ and  $S_1,S_2,\ldots,S_{\ell}$, respectively.
CHOLMOD determines the size of super-nodes by using 
heuristics like approximate minimum ordering so that 
the BLAS library can be used to process
the sparse Cholesky factorization and forward/backward substitution.
In addition, the structure of the memory allocated to primal
$\widehat{\L}$ is identical to that of dual $\N$
in MC-PDIPM.
Hence,  we first obtain the structure of $\N$ by constructing
the aggregate sparsity pattern $\AC$ and
applying CHOLMOD to obtain
its symbolic sparse Cholesky factorization,
then we extract its super-node information
($C_1,C_2,\ldots,C_{\ell}$ and $S_1,S_2,\ldots,S_{\ell}$) to prepare
the memory for $\widehat{\L}$.

Table~\ref{table:newcomp} shows the computation time reductions had by the
new factorization of $\widehat{\X}$.
The computation times of SDPA-C~6 and SDPA-C~7 were measured on
a relaxation of a max-clique problem on 
a lattice graph with the parameters $p=300$ and $q=10$.
The details of this SDP and the computation environment will be described
in Section~\ref{sec:results}.
The dimension of the
variable matrices $\X$ and $\Y$ was $n = p \times q = 3000$,
and the extended sparsity pattern $\EC$
was decomposed into 438 cliques, $C_1,\ldots, C_{438}$.
Since the maximum cardinality of the cliques,
$\max\{|C_r| : r=1,\ldots,438\}$, was only 59, 
the matrices were decomposed into small cliques.
We named three representative bottlenecks as follows:
S-ELEMENT is the time taken by
(\ref{eq:SCE}), (\ref{eq:newSCM}) or (\ref{eq:newSCM2})
to evaluate the SCM elements; 
S-CHOLESKY is the Cholesky factorization routine for the SCM,
and P-MATRIX is the computation of the primal auxiliary matrix
$\widehat{d\X}$ by (\ref{eq:dX}), (\ref{eq:newdX}) or (\ref{eq:newdX2}).
In addition, Sub-S-ELEMENTS and Sub-P-MATRIX are the times 
of the forward/backward substitutions 
like $\L^T \D \L \e_k$ or $\widehat{\L}^{-T} \widehat{\L}^{-1} \e_k$
in S-ELEMENTS and P-CHOLESKY. 

% \begin{table}[htbp]
% \begin{center}
% \caption{Computation time reduction due to the new factorization
% (the time unit is second.)}
% \label{table:newcomp}
% \begin{tabular}{|c|r|r|}
%  \hline
%  & SDPA-C~6.2.1 & SDPA-C~7.3.8 \\
%  \hline
%  S-ELEMENTS & 4205.70 & 2094.03 \\
%  S-CHOLESKY & 185.36 & 161.03 \\
%  P-MATRIX &  316.00 & 242.51 \\
%  Other & 22.22 & 17.93 \\
%  \hline
%  Total & 4729.28 & 2515.50 \\
%  \hline
%  \end{tabular}
%  \end{center}
% \end{table}
 
\begin{table}[htbp]
 \begin{center}
 \caption{Reduction in computation time due to the new factorization
(time in seconds)}
 \label{table:newcomp}
 \begin{tabular}{|l|r|r|}
  \hline
  & SDPA-C~6.2.1 & SDPA-C~7.3.8 \\
  \hline
  S-ELEMENTS & 4837.64 & 2054.98 \\
 \multicolumn{1}{|r|}{\hspace{0.5cm}(Sub-S-ELEMENTS)} & 3798.64 & 1512.22 \\
  \hline
  S-CHOLESKY & 167.65 & 161.23 \\
  \hline
  P-MATRIX &  346.46 & 242.06 \\
   \multicolumn{1}{|r|}{\hspace{0.5cm}(Sub-P-MATRIX)} & 322.54 & 224.16 \\
  \hline
  Other & 10.02 & 21.20 \\
  \hline
  Total & 5361.77 & 2479.49 \\
  \hline
  \end{tabular}
  \end{center}
 \end{table}

Table~\ref{table:newcomp} indicates that the new factorization reduced
the evaluation time of SCM (4837 seconds) to half (2054 seconds).
The computation time on P-MATRIX also shrank from 346 seconds to 242 seconds.
Consequently, the new factorization yielded a speedup of 2.16-times.
From the time reduction in Sub-S-ELEMENTS and Sub-P-MATRIX, 
%this is mainly because it does not require the
%diagonal-block matrix $\D$ as the original factorization and 
%it enables us to use the existing efficient package
%CHOLMOD. 
we can see that the removal of $\D$ enabled us to utilize the efficient 
forward/backward substitution of CHOLMOD.
As we will show in Section~\ref{sec:results}, 
the effect is even more pronounced for larger SDPs.

% At the end of this subsection, we note that when the overlaps among cliques $C_1,C_2,\ldots,C_{\ell}$
% are relatively large, some of cliques are merged into a clique for efficiency. For example,
% if $C_1 = \{1,2,\ldots,100\}, C_2 = \{2,3,\ldots,101\}, C_3 = \{3,4,\ldots,102\}$, these three cliques 
% are replaced by a new clique $C^{'}_1 = \{1,2,\ldots,101,102\}$. This merge is done with
% heuristic methods before the main loop of MC-PDIPM, and SDPA-C~7.3.8 uses the heuristics 
% implemented in CHOLMOD which matches $\widehat{\X}^{-1} = \widehat{\L} \widehat{\L}^T$, 
% while SDPA-C~6.2.1 uses a different heuristics which matches the factorization
% (\ref{eq:sparse_fact}).
% Hence,  even when we get $\widehat{\X} = \L\D\L^T$ in SDPA-C~6.2.1 and compute 
% $\widetilde{\L} = \L^{-T} \D^{-1/2}$
% with the matrix $\D^{1/2}$ such that $\D = \D^{1/2}\D^{1/2}$,
% we can not direclty put $\widetilde{\L}$ into the memory space of $\widehat{\L}$ of CHOLMOD.

\subsection{Matrix-completion primal-dual interior-point method
  for multithreaded parallel computing} \label{sec:threads}

To enhance the performance of SDPA-C~7 even further, we take
advantage of multithreaded parallel computing.
Since the processors on modern PCs have multiple cores (computation unit),
we can assign different threads (computation tasks) to the cores
and run multiple tasks simultaneously.
For example, multithreaded BLAS libraries are often used
to reduce the time related to linear algebra computations 
in a variety of numerical optimization problems.
However, the effect of multithreaded BLAS libraries is limited to 
dense linear algebra computations. Hence, we should seek a way to apply multithreaded
parallel computing 
to not only the dense linear algebra but also larger computation blocks of MC-PDIPM.

Here, we use multithreaded parallel computing to resolve the three bottlenecks of MC-PDIPM;
S-ELEMENTS, S-CHOLESKY, and P-MATRIX.
A parallel computation of these bottlenecks was 
already performed in SDPARA-C~\cite{NAKATA06}
with the MPI (Message Passing Interface) protocol on multiple PCs.
To apply multithreaded parallel computing on a single PC, however, 
we need a different parallel scheme.

% We start our discussion from S-ELEMENTS.
For S-ELEMENTS,
the SCM is evaluated with the formula~(\ref{eq:newSCM2}), 
and we can reuse the results of $\widehat{\L}^{-T}\widehat{\L}^{-1}\e_k$
and $\N^{-T}\N^{-1}\left[\A_j\right]_{*k}$ for
all the elements in $\left[\B\right]_{*j}$ 
(the $j$-th column of $\B$).
On the other hand, we can not reuse the results
for different columns of $\B$, 
since the memory needed to hold $\widehat{\L}^{-T}\widehat{\L}^{-1}\e_k$
for all $k=1,\ldots,n$ is equivalent to holding the fully dense matrix
$\widehat{\X}$, and this means we lose the nice sparse structure of $\EC$.

Since the computation of each column $\left[\B\right]_{*j}$
is independent from those of the other columns
$\{\left[\B\right]_{*i} : 1 \le i \le m, \quad i \ne j\}$,
SDPARA-C simply employs a column-wise distribution, in which the $p$th thread
evaluates the columns assigned by the set 
$\SC_p = \{j : 1 \le j \le m, \quad (j-1) \ {\tt mod}\  u = p - 1\}$,
where $a \ {\tt mod} \ b$ stands for the remainder of the division of $a$ by $b$ and
$u$ denotes the number of cores that 
participate in the multithreaded computation.

This simple assignment was necessary for SDPARA-C,
%as reported in \cite{NAKATA06},
since the memory taken up by $\B$ was assumed to be distributed 
over multiple PCs and we had to fix the column assignments in order not to send 
the evaluation result on one PC to another PC, which would have entailed a lot of network communications. 
In contrast, on a single PC, all of the threads of can share the memory space.
Hence, we devised more efficient parallel schemes
to improve the load-balance over all the threads.

\vspace{0.5cm}
\noindent {\bf Algorithm 2: Multithreaded parallel computing for the
evaluation of SCM}
\begin{enumerate}
 \item[Step 1] Initialize the SCM $\B=\O$.
	      Set $\SC = \{1,2,\ldots,m\}$.
 \item[Step 2] Generate $u$ threads.
 \item[Step 3] For $p=1,2,\ldots,u$, run the $p$th thread on the $p$th core
	      to execute the following steps:
	      \begin{enumerate}
	       \item[Step 3-1] If $\SC = \emptyset$, terminate.
	       \item[Step 3-2] Take the smallest element $j$ from
			    $\SC$ and update
			    $\SC \leftarrow \SC \backslash \{j\}$.
			    \item[Step 3-3] Evaluate
					 $\left[\B\right]_{*j}$ by \\
{\bf for} $k=1,\ldots, n$ \\
\hspace{0.5cm} Apply the forward/backward substitution routine
					 of CHOLMOD \\
\hspace{0.8cm}	to obtain
		$\v_1 := \widehat{\L}^{-T}\widehat{\L}^{-1} \e_k$ and
		$\v_2 := \N^{-T} \N^{-1}\left[\A_j\right]_{*k}$. \\
\hspace{0.5cm} {\bf for} $i=1,\ldots,m$ \\
\hspace{1.0cm} 	Update $B_{ij} \leftarrow B_{ij} + \v_1^T \A_i \v_2$.
	       \end{enumerate}
\end{enumerate}

\begin{figure}[htbp]
 \begin{center}
 \includegraphics[scale=0.7]{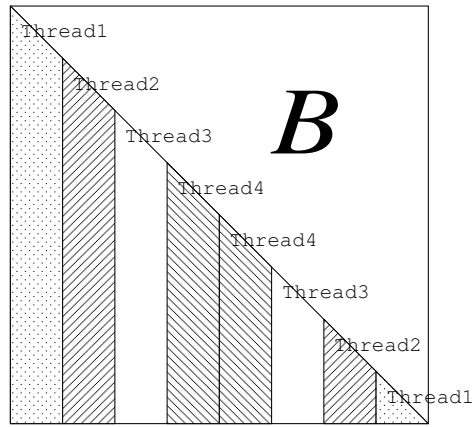}
 \caption{Thread assignment to the Schur complement matrix $\B$}
 \label{fig:Schur-Parallel}
  \end{center}
 \end{figure}

Figure~\ref{fig:Schur-Parallel} shows an example of thread
assignment to SCM $\B$ where $\B \in \SMAT^8$ and $u=4$.
Note that we evaluated only the lower triangular part of $\B$,
since $\B$ is symmetric.
We had $u=4$ threads; thus,
the $p$th thread evaluated the $p$th column for $p=1,2,3,\mbox{and}\ 4$
at the beginning.
Let us assume 
the computation cost is the same over all
$B_{ij}$ in Figure~\ref{fig:Schur-Parallel};
among the four threads, the 4th thread finishes its column evaluation
in the shortest time, so its evaluates the 5th column. 
After that, the 3rd thread finishes its
first task and moves to the 6th column.
On the other hand, when the 4th column requires a longer computation time than 
the 3rd column does, the 3rd thread takes the 5th column
and the 4th thread evaluates the 6th column.

% In this simple example, the number of elements $B_{ij}$
% evaluated by each thread is same.
% Consequently, this property generates a nice parallel performance of
% the new SDPA-C.

The overhead in Algorithm~2 is only in Step~3-2, where
only one thread should enter Step 3-2 at a time.
Hence, we expected that Algorithm~2
would balance the load better than the simple
column-wise distribution employed in SDPARA-C.
In particular, the main computation cost of each column
$\left[\B\right]_{*j}$ is
$\widehat{\L}^{-T}\widehat{\L}\e_k$ and
		$\N^{-T} \N^{-1}\left[\A_j\right]_{*k}$;
therefore, it is almost proportional to the number of 
nonzero columns of $\A_j$. When only a few of $\A_1, \ldots, \A_m$ 
have too-many nonzero columns and the others have only a few,
a simple column-wise distribution has a hard time keeping the
load-balanced. On the other hand, Algorithm~2 can naturally overcome this difficulty.

When we implement Algorithm~2,
we should pay attention to the number of
threads generated by the BLAS library that CHOLMOD internally calls for
the forward/backward substitution 
($\widehat{\L}^{-T}\widehat{\L}\e_k$ and
		$\N^{-T} \N^{-1}\left[A_j\right]_{*k}$).
For example, let us suppose that four cores are available ($u=4$).
If we generate four threads in Step~2 % of Algorithm~2 
and each thread internally
generates four threads for the BLAS library,
then we need to manage 16 threads in total on the four cores.
The overhead for this management is considerable,
and when we implemented the multithreaded parallel computing in
this way, SDPA-C took at least ten times longer 
than  single-thread computing. Therefore, we decided to turn off the
multithreading of the BLAS library before entering 
the forward/backward substitution routine
and turn it on again after the routine finishes.
%As the numerical result in Table~\ref{table:multiple-threads} will show later, 
%this scheme shrinks the computation time of the new SDPA-C~7.
% This enables us to resolve the thread conflicts.

Now let us examine S-CHOLESKY and P-MATRIX.
For S-CHOLESKY, our preliminary experiments indicated that the usage of the BLAS library 
for both LAPACK and MUMPS is sufficient for delivering the performance of multithreaded parallel computing.
In P-MATRIX, the primal auxiliary matrix $\widehat{d\X}$ was evaluated
by using formula~(\ref{eq:newdX2}).
Since this formula naturally indicates the independence of
$\widehat{d\X}$ columns, the simple column-wise distribution was 
employed in SDPARA-C. However,
in multithreading, all of the threads share the memory,
hence, we can replace the column-wise distribution
with the first-come first-served concept, the same parallel
concept as used in Algorithm~2.
We also used the above scheme to control the
number of threads involved in parallel computing.
% In addition, since $\widehat{d\X}$ is a fully dense matrix in general, each thread
% copies only the elements involved in 
% $\widehat{d\X}_{C_rC_r} \ (r=1,\ldots,\ell)$ after the evaluation of
% each column, so that we can compute the step length $\widehat{\alpha}_p$
% in (\ref{eq:primal_length2}) and update $\overline{\X}$.
% This etually enables us to execute the PD-IPM
% using the sparse matrices whose sparsity pattern is determined by $\EC$,
% % using the lower triangular matrix $\widehat{\L}$ that inherits the
% % sparsity of extended sparsity pattern 
% without completing $\overline{\X}$ to the fully dense matrix $\widehat{\X}$.

Table~\ref{table:multiple-threads} shows the computation time reduction
due to multithreaded parallel computing. It compares the results of SDPA-C~7
with those of  SDPA-C~6 and  SDPARA-C~1 on the same SDP.
% as Table~\ref{table:newcomp}.
In each bottleneck, the upper row is the computation time,
and the lower row is the speed-up ratio compared with a single thread.

%\begin{table}[htbp]
% \begin{center}
% \caption{Computation time reduction due to the multithreaded computing
%  (the time unit is second.)}
% \label{table:multiple-threads}
% \begin{tabular}{|l|r|r|r|r|r|r|r|}
%  \hline
%  & SDPA-C~6 &\multicolumn{3}{|c|}{SDPARA-C~1} & \multicolumn{3}{|c|}{SDPA-C~7}\\
%  \hline
%  threads & 1 & 1 & 2 & 4 & 1 & 2 & 4 \\
%  \hline
%  S-ELEMENTS & 4205.70 & 3103.51 & 1904.64 & 1309.88 & 2094.03 & 1170.37 & 731.98 \\
%             &         & 1.00x & 1.62x & 2.36x & 1.00x    & 1.78x & 2.86x \\
%  \hline
%  S-CHOLESKY & 185.36 & 312.41 & 110.76 & 65.12 & 161.03 & 85.24 & 47.77 \\
%             &        & 1.00x & 2.82x & 4.79x & 1.00x   & 1.88x  & 3.37x \\
%  \hline
%  P-MATRIX   & 316.00 & 316.04 & 159.32 & 81.34 & 242.51 & 131.03 & 83.54 \\
%             &        & 1.00x & 1.98x & 3.88x & 1.00x  & 1.85x & 2.90x \\
%  \hline
%  Other      & 22.22  & 32.80 & 44.76 & 40.69 & 17.93  & 23.88 & 25.86 \\
%  \hline
%  Total      & 4729.28 & 3764.84 & 2219.48 & 1497.03 & 2515.50 & 1410.52 & 889.15 \\
%             &         & 1.00x & 1.63x & 2.51x & 1.00x & 1.78x & 2.82x \\
%  \hline
% \end{tabular}
%\end{center}
%\end{table}
\begin{table}[htbp]
 \begin{center}
 \caption{Reduction in computation time  due to the multithreaded computing
  (time in seconds)}
 \label{table:multiple-threads}
 \begin{tabular}{|l|r|r|r|r|r|r|r|}
  \hline
  & SDPA-C~6 &\multicolumn{3}{|c|}{SDPARA-C~1} & \multicolumn{3}{|c|}{SDPA-C~7}\\
  \hline
  threads & 1 & 1 & 2 & 4 & 1 & 2 & 4 \\
  \hline
  S-ELEMENTS & 4873.64 & 3103.51 & 1904.64 & 1309.88 & 2054.98 & 1170.37 & 731.98 \\
             &         & 1.00x & 1.62x & 2.36x & 1.00x    & 1.76x & 2.81x \\
  \hline
  S-CHOLESKY & 167.65 & 312.41 & 110.76 & 65.12 & 161.23 & 85.24 & 47.77 \\
             &        & 1.00x & 2.82x & 4.79x & 1.00x   & 1.89x  & 3.38x \\
  \hline
  P-MATRIX   & 346.46 & 316.04 & 159.32 & 81.34 & 242.06 & 131.03 & 83.54 \\
             &        & 1.00x & 1.98x & 3.88x & 1.00x  & 1.85x & 2.90x \\
  \hline
  Other      & 10.02  & 32.80 & 44.76 & 40.69 & 21.06  & 23.88 & 25.86 \\
  \hline
  Total      & 5361.77 & 3764.84 & 2219.48 & 1497.03 & 2479.49 & 1410.52 & 889.15 \\
             &         & 1.00x & 1.63x & 2.51x & 1.00x & 1.76x & 2.79x \\
  \hline
 \end{tabular}
\end{center}
\end{table}

% In Table~\ref{table:multiple-threads},
The table shows that 
SDPA-C~7 with four threads reduced S-ELEMENTS to 731.98 seconds from 2054.98 seconds on a
single thread, a speedup of 2.81-times.
The computation time of P-MATRIX also shrank from 242.06 seconds
to 83.54 seconds, a speed up of 2.90-times.
These time reductions reduced the speed-up in the total time
from 2479.49 seconds to 889.15, for a speedup of 2.79-times.

SDPA-C~7 with 4 threads was 6.03-times faster than SDPA-C~6 using only a single thread.
The result indicate the parallel schemes
discussed above are effectively integrated into MC-PDIPM
and the new factorization $\widehat{\X}^{-1} = \widehat{\L}\widehat{\L}^T$.

The table also shows that SDPARA-C~1 took longer than SDPA-C~7.
This was mainly due to the overhead of the MPI protocol.
Since the MPI protocol is designed for multiple PCs, 
it is not appropriate for a single PC;
a multithreaded computation performs better.
% on a single PC. 
In addition, Algorithm~2 works more effectively in a multithreaded computing
environment than the simple column-wise distribution of SDPARA-C.
The speedup of SDPA-C~7 for S-ELEMENTS on four threads was
2.81, and it was higher than that of SDPARA-C~1 (2.36).

\section{Numerical Experiments}\label{sec:results}

We conducted a numerical evaluation of the performance of SDPA-C~7.
The computing environment was RedHat Linux run on a Xeon X5365 (3.0 GHz, 4 cores)
and 48 GB of memory space.
We used three groups of test problems, {\it i.e.},
max-clique problems over lattice graphs,
max-cut problems over lattice graphs, and spin-glass problems.
In this section, we will use the notation $\e$ to denote 
a vector of all ones, and $\e_i$ to denote a vector of all zeros except 1 at the $i$th element.

\vspace{0.5cm}
\noindent{\bf Max-clique problems over lattice graphs}

Consider a graph $G(V,E)$ with the vertex set $V = \{1,\ldots,n\}$
and the edge set $E \subset V \times V$.
A vertex subset $S \subset V$ is called a clique if
$(i,j) \in E$ for $\forall i \in S, \forall j \in S$.
The max-clique problem is to find a clique having the maximum cardinality
among all the cliques.% in $G(V,E)$.

Though the max-clique problem itself is NP-hard,
Lov\'{a}sz~\cite{LOVASZ79}
proposed an SDP relaxation method to obtain a good approximation 
in polynomial time.
The SDP problem below gives a good upper bound 
of the max-clique cardinality for $G(V,E)$,
\begin{eqnarray*}
 \max &:& \e\e^T \bullet \X \\
 \mbox{subject to} &:& \I \bullet \X = 1 \\
 & & \left(\e_i\e_j^T + \e_j\e_i^T\right) \bullet \X = 0
  \quad \mbox{for} \quad (i,j) \notin E \\
 & & \X \succeq \O.
\end{eqnarray*}

\begin{figure}[htbp]
 \begin{center}
 \includegraphics[scale=0.5]{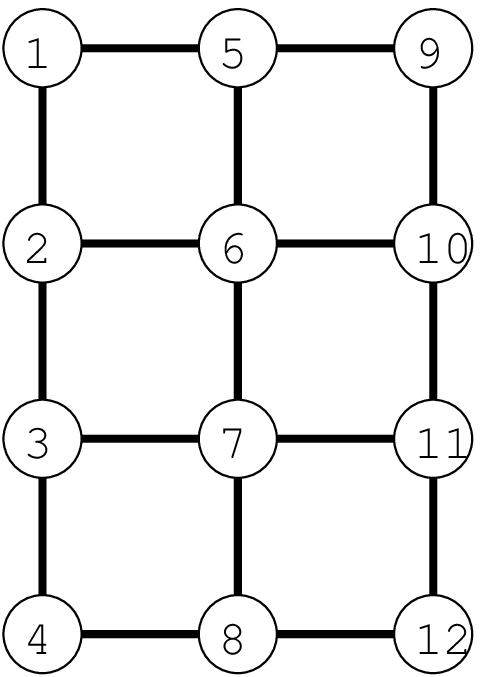}
 \caption{Lattice graph with the parameters $p = 4, q = 3$.}
 \label{fig:Lattice}
  \end{center}
 \end{figure}

For the numerical experiments, we generated SDPs of this type over
lattice graphs. A lattice graph $G(V,E)$ is determined by two
parameters $p$ and $q$, with the vertex set being
$V = \{1, 2, \ldots, p\times q\}$ and the edge set
$E = \{\left(i + (j-1) \times p , (i+1) + (j-1) \times p\right)
: i=1,\ldots,p-1, j=1,\ldots,q\}
\cup \{\left(i + (j-1) \times p, i + j \times p\right) : i=1,\ldots,p,
j=1,\ldots,q-1\}$.
An example of lattice graphs is shown in Figure~\ref{fig:Lattice},
where the parameters are $p = 4, q = 3$.

\begin{figure}[htbp]
 \begin{center}
\includegraphics[scale=0.4]{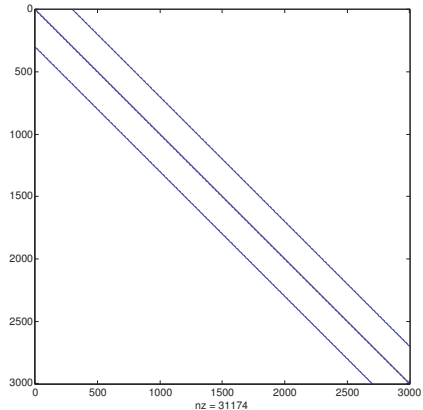}
 \caption{Aggregate sparsity pattern $\AC$
for the max-clique SDP with $p = 300, q = 10$.}
 \label{fig:max-clique}
  \end{center}
 \end{figure}

We applied the pre-processing technique proposed in \cite{FUKUDA00}
to convert the above SDP into an equivalent but sparser SDP.
Figure~\ref{fig:max-clique} shows the aggregate sparsity pattern $\AC$ 
for the max-clique SDP with $p=300, q=10$.
We applied approximate minimum degree heuristics to $\AC$ to make
Figure~\ref{fig:max-clique},  and  this figure
shows the sparse structure embedded in this SDP.
The sizes of cliques $C_1,\ldots,C_{\ell}$ 
can be much smaller in comparison with $n = p \times q = 300 \times 10 = 3000$, and
this $\AC$ does not incur
any fill-in, that is, $\EC = \AC$.
This SDP was the example solved in Section~\ref{sec:new}.

\vspace{0.5cm}
\noindent{\bf Max-cut problems over lattice graphs}

The SDP relaxation method for solving max-cut problems due to Goemans and
Williamson~\cite{GOEMANS95} is well-known, and it marked the beginning of  studies on
SDP relaxation methods.
Here, one considers a graph $G(V,E)$ 
with the vertex set $V = \{1,\ldots,n\}$
and the edge set $E \subset V \times V$.
Each edge $(i,j) \in E$ has the corresponding non-negative
weight $w_{ij}$ (for simplicity, $w_{ij} = 0$ if $(i,j) \notin E$).
The weight of the cut $C \subset V$ is the total weight of the edges 
traversed between $C$ and $V \backslash C$.
The max-cut problem is to find a subset $C$ which maximizes the cut
weight,
\begin{eqnarray*}
 \max&:& \sum_{i \in C, j \in V \backslash C} w_{ij} \qquad 
 \mbox{subject to} \ :\ C \subset V.
\end{eqnarray*}
An SDP relaxation of this problem is given by
\begin{eqnarray*}
 \min &:&  \A_0 \bullet \X \\
 \mbox{subject to} &:& \left(\e_i \e_i^T\right) \bullet \X = 1
  \quad \mbox{for} \quad i=1,\ldots,n \\
 & & \X \succeq \O,
\end{eqnarray*}
where $\A_0$ is defined by $ A_0 = \left(-\mbox{diag}(\W \e) +\W\right)$,
$\W$ is the matrix whose $(i,j)$ element is $w_{ij}$ for
$i=1,\ldots,n,\ j=1,\ldots,n$, 
% $\e$ is the vector in which all elements are 1,
and $\mbox{diag}(\w)$ is the diagonal matrix whose diagonal elements are 
the elements of $\w$.

When we generate SDP problems from the max-cut problem over the lattice
graphs, its aggregate sparsity pattern appears in 
the coefficient matrix of the objective function $\A_0$.
Hence, we can find a similar structure to the one shown in  Figure~\ref{fig:max-clique}
in its aggregate sparsity pattern.

\vspace{0.5cm}
\noindent{\bf Spin-glass problems}

The four SDPs of this type were collected  as a torus set
in the 7th DIMACS benchmark problems~\cite{DIMACS}.
These SDPs arise in computations of the ground-state energy 
of Ising spin glasses in quantum chemistry.
More information on this energy computation can be found at the 
Spin Glass Server~\cite{SPINGLASS} webpage and references therein. 

\begin{figure}[htbp]
 \begin{center}
\includegraphics[scale=0.4]{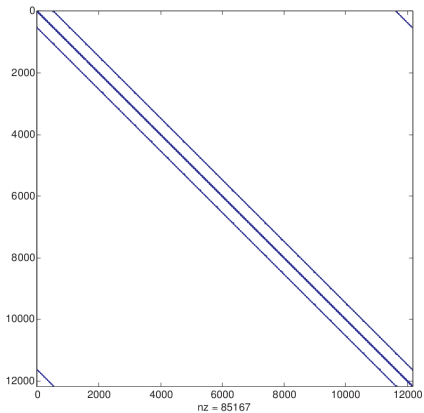}
 \caption{The aggregate sparsity pattern $\AC$ 
of the spin-glass SDP with $p = 23$.}
 \label{fig:spin-glass}
  \end{center}
 \end{figure}

The Ising spin-glass model has a parameter $p$
(the number of samples), 
and if we generate an SDP from a 3D spin-glass model, 
the dimension of the variable
matrices $\X$ and $\Y$ is $n = p^3$ \cite{SPINGLASS}.
Figure~\ref{fig:spin-glass} illustrates the aggregate sparsity pattern
$\AC$ of the spin-glass SDP with $p=23$ and $n = 23^3 = 12167$.
% We can again find a nice structure that matches the MC-PDIPM.

\vspace{0.5cm}
Table~\ref{table:SDPs} summarizes the SDPs of
the numerical experiments.
The first column is SDP's name and 
the second $p$ is the parameter used to generate it
(We fixed the parameter $q$ to $10$ for the max-clique problems and
the max-cut problems).
The third column $n$ is the dimension of the variable matrices
$\X$ and $\Y$, and
the fourth column
is the density of aggregate sparsity
pattern defined by $\frac{|\AC|}{n^2}$.
The fifth column $\ell$ is the number of cliques ($C_1,\ldots,C_{\ell}$),
and the sixth and seventh columns are the average and maximum sizes 
of the cliques defined by $\frac{\sum_{r=1}^{\ell}|C_r|}{\ell}$ and
$\max_{r=1,\ldots,\ell}C_r$, respectively.
% As described in Section~\ref{sec:decomp}, we extract the clique
% structure $C_1,\ldots,C_{\ell}$ from the super-node information of CHOLMOD.
The eighth column $m$
is the number of input data matrices $\A_1, \ldots, \A_m$.

\begin{table}[htbp]
 \begin{center}
  \caption{The sizes of the SDPs in the numerical experiments}
 \label{table:SDPs}
 \begin{tabular}{|l|r|r|r|r|r|r|r|r|r|r|}
  \hline
  Name & \multicolumn{1}{|c|}{$p$} & \multicolumn{1}{|c|}{$n$} & density  & \multicolumn{1}{|c|}{$\ell$} & ave-size & max-size  & \multicolumn{1}{|c|}{$m$} \\
  \hline
  MaxClique300 & 300 & 3000  & 0.38\%  & 348 & 28.36 & 51 & 5691 \\
  MaxClique400 & 400 & 4000  & 0.28\% & 439 & 29.89 & 59 & 7591\\
  MaxClique500 & 500 & 5000  & 0.23\% & 581 & 28.26 & 50 & 9491\\
  \hline
  MaxCut400 & 400 & 4000  & 0.15\% &  1282 & 8.05 & 26 & 4000 \\
  MaxCut500 & 500 & 5000  & 0.12\% & 1607 & 8.04 & 26 & 5000 \\
  MaxCut600 & 600 & 6000  & 0.097\% & 1932 & 8.04 & 26 & 6000 \\
  MaxCut800 & 800 & 8000  & 0.072\% & 2582 & 8.03 & 26 & 8000 \\
  MaxCut1000 & 1000 & 10000  & 0.058\% & 3232 & 8.02 & 26 & 10000 \\
  MaxCut1200 & 1200 & 12000  & 0.048\% & 3882 &  8.02 & 26 & 12000 \\
  \hline
  SpinGlass10 & 10 & 1000  & 0.80\% & 155 & 25.69 & 294 & 1000\\
  SpinGlass15 & 15 & 3375  & 0.24\% & 191 & 29.97 & 773 & 3375\\
  SpinGlass18 & 18 & 5832  & 0.13\% & 1118 & 28.13 & 913 & 5832\\
  SpinGlass20 & 20 & 8000  & 0.10\% & 1737 & 25.75 & 1080 & 8000\\
  SpinGlass23 & 23 & 12167  & 0.066\% & 2556 & 27.30 & 1488 & 12167\\
  SpinGlass25 & 25 & 15625  & 0.051\%  & 3173 & 29.50 & 1798 & 15625\\
  \hline
\end{tabular}
\end{center}
\end{table}

% In this section, 
We compared the computation times 
of SDPA-C~6~\cite{NAKATA06}, SDPA-C~7, and
SDPA~7~\cite{SDP-HANDBOOK} and SeDuMi~1.3~\cite{STURM99}.
The former two implemented MC-PDIPM, while the latter two implemented
the standard PDIPM.
Here, we did not conduct a numerical experiment on SDPARA-C, since we found that 
the overhead due to the MPI protocol 
was a severe disadvantage when we ran it on a single PC
as shown in Section~\ref{sec:threads}.

Table~\ref{table:SDP-time} lists
the computation times of the four solvers using their default parameters.
We used four threads for SDPA-C~7 and SDPA~7.
% We also assign the 4 threads to the BLAS libarary, hence,
% SDPA-C6 and SeDuMi use the 4 threads when they internally call the
% BLAS library.
The symbol '$>$2days' in the table indicates that 
we gave up on the SeDuMi execution since it required at least two days.

\begin{table}[htbp]
 \begin{center}
  \caption{Computation times of four solvers for the SDPs
  in Table~\ref{table:SDPs} (time in seconds)}
 \label{table:SDP-time}
 \begin{tabular}{|l|r|r|r|r|}
  \hline
  Name & SDPA-C6 & SDPA-C7 & SDPA~7 & SeDuMi1.3 \\
  \hline
  MaxClique300 & 4792.28  & 889.15  & 11680.12  & 10260.15 \\
  MaxClique400 & 12681.16 & 1903.13 & 26159.40 & 24824.05 \\
  MaxClique500 & 19973.98 & 3733.41 & 38265.02  & 46168.56 \\
  \hline
  MaxCut400 & 386.80 & 539.22 & 3686.78 &  16779.68\\
  MaxCut500 & 683.27 & 876.81 & 6548.26 & 32557.13 \\
  MaxCut600 & 1194.89 & 1295.01 & 11098.35 & 60444.33 \\
  MaxCut800 & 2518.80 & 2371.43 & 25377.62 & 146235.81 \\
  MaxCut1000 & 4301.11 & 4032.80 & 47270.45 & $>$2days \\
  MaxCut1200 & 7400.28 & 6030.00 & 75888.85 & $>$2days \\
  \hline
  SpinGlass10 & 50.10 & 20.77 & 11.85 & 228.71 \\
  SpinGlass15 & 1306.00 & 560.40 & 336.40 & 13789.58 \\
  SpinGlass18 & 6734.40 & 2136.88 & 1522.60 & 68570.89 \\
  SpinGlass20 & 15450.13 & 4552.10 & 3726.03  & $>$2days \\
  SpinGlass23 & 55942.57 & 13184.25  & 12598.14  & $>$2days \\
  SpinGlass25 & 107502.48 & 24913.20 & 26023.67 & $>$2days \\
  \hline
\end{tabular}
\end{center}
\end{table}

% We now examine the details of Table~\ref{table:SDP-time}.
On the max-clique problems, the MC-PDIPM solvers were faster than
the standard PDIPM solvers. 
Since the matrix-completion method benefited from  the nice properties of lattice graphs,
even SDPA-C~6 was twice as faster as SDPA~7.
The detailed breakdown of the time on MaxClique400 is displayed in Table~\ref{table:mcq400}
(since SeDuMi does not print out its internal computation time,
we did not list its breakdown).
As shown in the SDPA~7 column, the standard PDIPM took a long time on P-MATRIX (\ref{eq:dX})
and Other (mainly, the computation of the step length by (\ref{eq:primal_length})).
Though these parts required an $O(n^3)$ computation cost,
the MC-PDIPM decomposed the full matrix $\X$ into the sub-matrices
$\overline{\X}_{C_rC_r} \ (r=1,\ldots,\ell)$;
hence, it was able to reduce the computation cost of these two parts.
Furthermore, SDPA-C~7 resolved the heaviest parts of SDPA-C~6
by using the new factorization
$\widehat{\X}^{-1} = \widehat{\L}\widehat{\L}^T$ and
multithreaded parallel computing.
% as discussed in Section~\ref{sec:new}.
Consequently, SDPA-C~7 was the fastest among the four solvers;
in particular, it was 12.36-times faster
than SeDuMi on MaxClique500. 

\begin{table}[htbp]
 \begin{center}
 \caption{Computation times on MaxClique400
(time in seconds)}
 \label{table:mcq400}
 \begin{tabular}{|c|r|r|r|r|}
  \hline
  & SDPA-C~6 & SDPA-C~7 & SDPA~7 & SeDuMi~1.3 \\
  \hline
  S-ELEMENTS & 9314.34 & 1431.69 & 262.04 & --\\
  S-CHOLESKY & 2601.04 & 248.55 & 403.28 & -- \\
  P-MATRIX & 734.41 & 175.24 & 13151.10 & -- \\
  Other & 31.37 & 47.65 & 12342.98 & -- \\
  \hline
  Total & 12681.16 & 1903.13 & 26159.40 & 24824.05 \\
  \hline
  \end{tabular}
  \end{center}
 \end{table}

For the max-cut problems,
though MC-PDIPM was again superior to the standard PDIPM,
SDPA-C~7 was less effective than SDPA-C~6.
In particular, it took longer on S-ELEMENTS and
P-MATRIX, both of which utilized multithreaded computing.
It needed an overhead to generate the threads, and the input
matrices of the max-cut problems were too simple
to derive any benefit from multithreaded computing.
Indeed, each input matrix $\A_i = \e_i\e_i^T$ has only one nonzero element,
and this is reflected in the short computation time of SDPA~7's S-ELEMENTS. 
In the standard PDIPM,
the S-ELEMENTS computation is an inexpensive task of (\ref{eq:SCM})
with $\A_i = \e_i\e_i^T$ and $\A_j = \e_j\e_j^T$,
since the fully dense matrices $\X$ and $\Y^{-1}$ are obtained
with extensive memory space and a heavy computation through
P-MATRIX and the inverse of the fully dense matrix.
For the large max-cut problems, however,  SDPA-C~7 solved the
SDPs faster than SDPA-C~6 or SDPA~7. As shown in the Max1200 result of 
Table~\ref{table:mc500}, SDPA-C~7 still incurred a multithreading overhead on
 S-ELEMENTS and P-MATRIX, but the multithreaded 
BLAS library resolved the principal bottleneck, S-CHOLESKY.
We can say that SDPA-C~7 would work even better on larger SDPs of this type.

\begin{table}[htbp]
 \begin{center}
 \caption{Computation times on MaxCut500 and MaxCut1200
(time in seconds)}
 \label{table:mc500}
 \begin{tabular}{|c|r|r|r|r|}
  \hline
  \multicolumn{5}{|c|}{MaxCut500} \\
  \hline
  & SDPA-C~6 & SDPA-C~7 & SDPA~7 & SeDuMi~1.3 \\
  \hline
  S-ELEMENTS & 105.25 & 317.56 & 16.17  & --\\
  S-CHOLESKY & 502.43 & 226.06 & 214.43  & -- \\
  P-MATRIX & 61.21 & 315.73 & 3254.05  & -- \\
  Other & 14.38 & 17.46 & 3063.60  & -- \\
  \hline
  Total & 683.27 & 876.81 & 6548.26 &  32557.13 \\
  \hline
  \multicolumn{5}{|c|}{MaxCut1200} \\
  \hline
  & SDPA-C~6 & SDPA-C~7 & SDPA~7 & SeDuMi~1.3 \\
  \hline
  S-ELEMENTS & 1265.02 & 1800.27 & 107.17  & --\\
  S-CHOLESKY & 5562.48 & 2318.17 & 2674.50  & -- \\
  P-MATRIX & 532.21 & 1837.53 & 39490.81  & -- \\
  Other & 40.57 & 74.03 & 33616.37  & -- \\
  \hline
  Total & 7400.28 & 6030.00 & 75888.85 & $>$2days \\
  \hline
  \end{tabular}
  \end{center}
 \end{table}

SDPA~7 was the fastest  in solving the spin-glass SDPs with $p=10$,
since its standard PDIPM is more effective on smaller SDPs 
where the variable matrices are small and do not need to be decomposed.
When we increased $p$, however, the time difference
between SDPA-C~7 and SDPA~7  shrank, and SDPA-C~7 became faster than SDPA~7  at $p=25$.
The reason why the growth in the computation time of SDPA-C~7 was not as steep 
in comparison with SDPA~7 is that the average size of cliques does not grow 
with $p$, as shown in Table~\ref{table:SDPs}.
In particular, as seen in the breakdown of the computation times in
Table~\ref{table:spinglass18},
this affects P-MATRIX % (the computation of $\Delta \X$ in (\ref{eq:newdX})),
and its computation time in MC-PDIPM grows more gradually 
than in the standard PDIPM.
SpinGlass25 was the largest among the spin-glass SDPs
in our experiments; it required
almost 48~GB of memory, 
close to the capacity of our computing environment.
We expect that SDPA-C~7 would be more effective on larger SDPs 
of the spin-glass type.
The ratios of SpinGlass25 over SpinGlass23
were $\frac{107502.48}{55942.57} = 1.92$ in SDPA-C~6, 
$\frac{24913.20}{13184.25} = 1.89$ in SDPA-C~7,
and $\frac{26023.67}{12598.14} = 2.07$ in SDPA~7.
In addition, SDPA~7 incurred a considerable computation cost on the other parts,
{\it i.e.}, 'Other'. This 'Other' category in SDPA~7 contained
the miscellaneous parts related to
the computation of the variable matrices $\X$ and $\Y$.
Since they are miscellaneous and we can not say which costs the most,
we did not examine them in detail,
but we note that the fully dense properties of $\X$ and $\Y$
diminished the performance of the 'Other' parts of SDPA~7.

We should emphasize that MC-PDIPM is not the only reason for 
SDPA-C~7 being faster, because Table~\ref{table:spinglass18} shows
that SDPA-C~6 was much slower than SDPA~7.
The new factorization of $\widehat{\X}^{-1} =
\widehat{\L}\widehat{\L}^T$ 
and the multithreaded computing were the keys to solving the spin-glass
SDPs in the shortest time.
% Furthermore, the employiment of CHOLMOD improves the above ratio from
% 1.92 in SDPA-C~6 to 1.89 in SDPA-C~7.

\begin{table}[htbp]
 \begin{center}
 \caption{Computation times on Spinglass18 and Spinglass25
(time in seconds)}
 \label{table:spinglass18}
 \begin{tabular}{|c|r|r|r|r|}
  \hline
  \multicolumn{5}{|c|}{Spinglass18} \\
  \hline
  & SDPA-C~6 & SDPA-C~7 & SDPA~7 & SeDuMi~1.3 \\
  \hline
  S-ELEMENTS & 2738.18 & 1012.13 & 22.93  & --\\
  S-CHOLESKY & 178.15 & 52.70 & 30.69  & -- \\
  P-MATRIX & 2516.52 & 984.76 & 620.54  & -- \\
  Other &  1301.5& 87.29 & 848.44  & -- \\
  \hline
  Total & 6734.40 & 2136.88 & 1522.60 &  68570.89 \\
  \hline
  \multicolumn{5}{|c|}{Spinglass25} \\
  \hline
  & SDPA-C~6 & SDPA-C~7 & SDPA~7 & SeDuMi~1.3 \\
  \hline
  S-ELEMENTS & 35314.96 & 11829.56 & 220.19  & --\\
  S-CHOLESKY & 3681.44 & 945.49 & 520.39  & -- \\
  P-MATRIX & 45238.37 & 11594.76 & 11846.59  & -- \\
  Other &  23267.71 &  588.39 & 13436.50  & -- \\
  \hline
  Total & 107502.48 & 24913.20 & 26023.67 & $>$2days \\
  \hline
  \end{tabular}
  \end{center}
 \end{table}

%ここで，Spinglass18 で SDPA~7 の Other の内訳は
%$\mu \I - \X \Y$ が 320.51 秒, $\alpha$ の計算が 145.35 秒，
%$\X,\Z$ のコレスキー分解がそれぞれ 117.07, 198.73 秒．
%これらを Other に含んでいるので，そのあたりをどう表現するか．
%また，P-MATRIX は SDPA-C~7 の方が SDPA-7 よりもゆっくりと
%上昇していて，このあたりも文章に書いてもいいと思われる．

Finally, Table~\ref{table:SDP-memory} shows the amount of memory required
to solve the SDPs in Tables~\ref{table:mcq400}, \ref{table:mc500}, and
\ref{table:spinglass18}.
The notation '$>$ 31G' indicates that SeDuMi exceeded
the time limit (2 days) and used 31~gigabytes of memory
during the two-day execution.
By comparison, MC-PDIPM saved a lot of memory
by removing the fully dense matrices.
For example, in MaxClique400, SDPA-C7 used
only $\frac{1}{6}$ times and $\frac{1}{10}$ times 
the memory  of SDPA7 and SeDuMi, respectively.
In addition, the new factorization
reduced the memory needed for the largest SDP (Spinglass25) from
8.1 gigabytes in SDPA-C6 to 3.7 gigabytes in SDPA-C7.
It reduced the required memory because it can reuse the memory structure of CHOLMOD.
% it is effective to reduce the required memory space.

\begin{table}[htbp]
 \begin{center}
  \caption{Amount of memory space required for solving SDPs
  in Tables~\ref{table:mcq400}, \ref{table:mc500}, and \ref{table:spinglass18} (M and G indicate megabytes and gigabytes.)}
 \label{table:SDP-memory}
 \begin{tabular}{|l|r|r|r|r|}
  \hline
  Name & SDPA-C6 & SDPA-C7 & SDPA~7 & SeDuMi1.3 \\
  \hline
  MaxClique400 & 548M & 516M & 3.0G  &  5.2G \\
  \hline
  MaxCut500 & 249M & 236M & 4.1G &  6.1G  \\
  MaxCut1200 & 1.2G & 1.2G & 23G & $>$ 31G  \\
  \hline
  SpinGlass18 & 2.0G & 707M & 5.6G & 8.3G \\
  SpinGlass25 & 8.1G & 3.7G & 40G & $>$ 36G \\
  \hline
\end{tabular}
\end{center}
\end{table}

\section{Conclusions and Future Directions}\label{sec:conclusions}

We implemented a new SDPA-C~7, that uses 
a more effective factorization of $\widehat{\X}^{-1} =
\widehat{\L}\widehat{\L}^T$ 
and takes advantages of multithreaded parallel computing.
Our numerical experiments verified that these two improvements
enhanced the performance of MC-PDIPM
and reduced  the computation times of
the max-clique and spin-glass SDPs.

SDPA-C~7 is available at the SDPA web site,
http://sdpa.sourceforge.net/.
Unlike SDPA-C~6, SDPA-C~7 has a callable library
and a Matlab interface;  it can now be embedded in other C++
software packages and be directly called from inside Matlab.
The callable library and Matlab interface
will no doubt expand the usage of SDPA-C.

%Finally, we discuss some future directions.
As shown in the numerical experiments on the max-cut problems,
if the input SDP has a very simple structure, we should automatically 
turn off the multithreading. However, this would require
a complex task to estimate the computation time accurately over the
multiple threads from the input SDPs.
Another point is that SDPA-C~7 has a tendency to be faster for large
SDPs. This is an excellent feature,  
but it does not extend to smaller SDPs. 
Although this is mainly because MC-PDIPM is 
intended to solve large SDPs with the factorization of the variable matrices,
% based on the structural sparsity,
we should combine it with other methods that effectively compute 
the forward/backward substitution 
% of  sparse matrices 
of small dimensions.

\section*{Acknowledgments}
The authors thank Professor Michael J\"{u}nger of Universit\"{a}t
zu K\"{o}ln and 
Professor Frauke Liers of Friedrich-Alexander Universit\"{a}t
Erlangen-N\"{u}rnberg for providing us with the general instance generator
of the Ising spin glasses computation.
The authors gratefully acknowledge the constructive comments of the 
anonymous referee.

\end{document}